\documentclass[12pt]{amsart}

\newcommand{\Vect}{\mathop{\rm Vect\,}\nolimits}
\newcommand{\Spec}{\mathop{\rm Spec}\nolimits}
\newcommand{\Gal}{\mathop{\rm Gal\,}\nolimits}
\newcommand{\tr}{\mathop{\rm tr\,}\nolimits}
\newcommand{\Z}{\mathbf Z}
\newcommand{\C}{\mathbf C}
\newcommand{\F}{\mathbf F}
\newcommand{\T}{\mathbf T}
\newcommand{\G}{\mathbf G}

\renewcommand{\P}{{\rm P}}

\newcommand{\GL}{{\rm GL}}
\newcommand{\Gr}{{\rm Gr}}
\newcommand{\Fl}{{\rm Fl}}

\begin{document}
\title{Some problems in mathematics and mathematical physics}
\author{Alexander Roi Stoyanovsky}
\begin{abstract}
We discuss new approaches to fundamental problems of mathematics and mathematical physics such as mathematical foundation of quantum field theory, the Riemann hypothesis, and construction of
noncommutative algebraic geometry.
\end{abstract}
\maketitle

The purpose of this paper is to state some mathematical problems related to mathematical foundation of quantum field theory, to the Riemann hypothesis on zeros of the Riemann zeta function, and to 
constructing noncommutative algebraic geometry.  

\section{Mathematical quantum gravity}

This problem arose from 20 year attempts to answer the following question: what is mathematics behind quantum field theory?

A naive approach to this question is to try to generalize the wave theory (the theory of linear partial differential equations) to the case of multidimensional variational principles and multidimensional bicharacteristics. 
The result of author's attempts in this direction is the conjecture [5] that perturbative and axiomatic quantum field theory is not a mathematical theory, i.~e. that axioms of Lorentz covariant quantum field
theory cannot be mathematically realized, like the Continuum Hypothesis. In other words, mathematical synthesis of quantum theory and special relativity does not exist.   

However, there is an example of mathematics behind quantum field theory, named two-dimensional conformal field theory (see e.~g. [1, 3, 4] and references therein). 
2D conformal field theory is a quantum field theory for the spacetime being a complex curve, invariant with respect to conformal transformations, i.~e. (complex analytical) diffeomorphisms of the curve, or at least with 
respect to the Lie algebra of (holomorphic) vector fields on the curve.  

This leads to an assumption that, while quantum field theory with ``small'' symmetry groups cannot be mathematically founded, quantum field theory with ``large'' symmetry groups can be stated mathematically.

The question is: what is the analog of 2D conformal field theory for the case of 
smooth real spacetime manifold $M$ of arbitrary dimension?  

It should be a quantum field theory invariant with respect to the group of diffeomorphisms of $M$, 
or at least with respect to the Lie algebra $\Vect M$ of (smooth) vector fields on $M$. 

An example of classical field theory invariant with respect to the group of diffeomorphisms of $M$ is the Einstein general relativity. 
Thus, the question is to construct mathematical quantization of general relativity, i.~e. {\it mathematical quantum gravity}.  

Quantum field theory is usually defined as an algebra of observables with operator product expansion, and the space of states on which this algebra acts, see e.~g. [1,~2]. 
Note that the space of states in [2] is not a vector space. We do not know whether there should exist
appropriate space of states in mathematical quantum gravity. We expect that, in general, it does not exist, unlike in topological field theory.
{\it Our problem is}: to construct only a $\Vect M$-invariant deformation quantization of the Poisson algebra of local observables in general relativity, with 
the structure of $\Vect M$-invariant operator product expansion.

\section{Hypergeometric zeta functions}

Here we present a new approach to fundamental problems of algebraic geometry and number theory such as the Riemann hypothesis.

It is known that modern ``abstract'' methods do not work for the Riemann hypothesis. The reason for this seems to be very simple:  modern methods  
use only algebraic operations $+,-,\times,:$, and they do not use the analytic operations like derivative and integral. Some modern deep methods in algebraic geometry, like Hodge theory or the theory of $D$-modules and constructible 
sheaves, use differential equations, but they use only existence of solutions, and do not use the essential quantitative properties of solutions. The main problem of modern algebraic geometry is that it has poor amount of functions. 
In order to prove the Weil conjectures, one invented \'etale cohomology of schemes, whose main idea is to enrich the amount of algebraic functions and include monodromy and constructible sheaves. However, for 
fundamental problems such as the Riemann hypothesis this amount of functions is still too small.  

In a large series of papers (see in particular [6] and references therein), I.~M.~Gel\-fand et al. invented and studied the astonishing rich world of general hypergeometric functions. 
The main idea of our approach to algebraic geometry and number theory is to revise these disciplines in terms of (appropriately developed) theory of general hypergeometric functions. We hope that 
hypergeometric functions can do the same work and are sometimes even better than abstract cohomological methods.

For example, one of the main problems of modern algebraic methods is that there is no general algebraic formula for the root of an algebraic equation
\begin{equation}
P(x)\equiv c_m x^m+c_{m-1}x^{m-1}+\ldots+c_0=0.
\end{equation}
In algebra one adjoins roots formally and studies their algebraic properties. However, there is an analytic formula for the root:
\begin{equation} 
x=\frac1{2\pi i}\oint\frac{P'(t)}{P(t)}tdt=-\frac1{2\pi i}\oint\log P(t)dt.
\end{equation} 

This formula implies that {\it the root of algebraic equation \emph{(1)} is a \emph(multi-valued\emph)  hypergeometric function of the coefficients} $c_0$, $\ldots$, $c_m$. 
This fact (in a slightly different form) has been discovered by Sturmfels [7].  It is the ``constructive fundamental theorem of algebra''.

As another illustration, we shall give a generalization of zeta functions of algebraic varieties over finite fields, based on the notion of hypergeometric function over a finite field defined in [6]. 
We call these generalizations by hypergeometric zeta functions. 

Let 
\begin{equation}
P(t_1,\ldots,t_n)=\sum\limits_{\omega=(\omega_1,\ldots,\omega_n)\in\Z^n}c_\omega t_1^{\omega_1}\ldots t_n^{\omega_n}
\end{equation}
be a Laurent polynomial in $n$ variables with coefficients $c_\omega$ from a finite field $\F_q$ with $q$ elements. Let $\chi:\F_q\to\C^\times$ be a fixed nontrivial additive character of $\F_q$, and let
$\pi_1,\ldots,\pi_n:\F_q^\times\to\C^\times$ be multiplicative characters, where $\F_q^\times=\F_q\setminus\{0\}$, $\C^\times=\C\setminus\{0\}$. The {\it general hypergeometric function over $\F_q$}
has been defined in [6] as
\begin{equation}
\Phi_q(\chi;P;\pi_1,\ldots,\pi_n)=\sum\limits_{t_1,\ldots,t_n\in\F_q^\times}\chi(P(t_1,\ldots,t_n))\pi_1(t_1)\ldots\pi_n(t_n).
\end{equation}

Let  $P_1,\ldots,P_k$ be Laurent polynomials in $t_1,\ldots,t_n$, let $\lambda_1,\ldots,\lambda_k$ be additional variables, let 
\begin{equation}
P(\lambda_1,\ldots,\lambda_k,t_1,\ldots,t_n)=\lambda_1P_1(t_1,\ldots,t_n)+\ldots+\lambda_kP_k(t_1,\ldots,t_n)
\end{equation}
(the {\it Cayley trick}), and let 
\begin{equation}
\pi_1=\ldots=\pi_{k+n}\equiv1
\end{equation}
 be the trivial characters. Denote 
\begin{equation}
\Phi_q(\chi;P;\pi_1,\ldots,\pi_{k+n})=\Psi_q(\chi;P_1,\ldots,P_k).
\end{equation}
Then we have 
\begin{equation}
\begin{aligned}
&\sum\limits_{\begin{subarray}{c}\lambda_1,\ldots,\lambda_k\in\F_q \\ t_1,\ldots,t_n\in\F_q^\times\end{subarray}}\chi(\lambda_1P_1(t_1,\ldots,t_n)+\ldots+\lambda_k P_k(t_1,\ldots,t_n))\\
&=(q-1)^n+\sum\limits_{\begin{subarray}{c}0<l\le k\\ 1\le i_1<\ldots<i_l\le k\end{subarray}}\Psi_q(\chi;P_{i_1},\ldots,P_{i_l})=q^k\#X(\F_q),
\end{aligned}
\end{equation}
where $\#X(\F_q)$ is the number of elements in the set $X(\F_q)$ of $\F_q$-points of the algebraic variety $X$ over $\F_q$ given by the system of equations and inequalities
\begin{equation}
t_1\ne0,\ldots,t_n\ne0,\ \ P_1(t_1,\ldots,t_n)=\ldots=P_k(t_1,\ldots,t_n)=0,
\end{equation}
i.~e. it is the number of solutions $(t_1,\ldots,t_n)\in\F_q^n$ of system (9).

Recall that the zeta function of an algebraic variety $X$ over $\F_q$ is defined as
\begin{equation}
\zeta_X(s)=\prod\limits_{x\in X}(1-\#k(x)^{-s})^{-1}=\exp\sum\limits_{d=1}^\infty\frac{q^{-ds}}d\#X(\F_{q^d}),
\end{equation}
where $\#k(x)=q^{d(x)}$ is the order of the residue field $k(x)=\F_{q^{d(x)}}$ of the point 
\begin{equation}
x\in X=\bigcup\limits_{d=1}^\infty X(\F_{q^d})/\Gal(\F_{q^d}/\F_q);
\end{equation}
here $\Gal(\F_{q^d}/\F_q)$ is the Galois group of $\F_{q^d}$ over $\F_q$.

Define the {\it hypergeometric zeta function} corresponding to a Laurent polynomial $P(t_1,\ldots,t_n)$ and to characters $\chi,\pi_1,\ldots,\pi_n$, as
\begin{equation}
\zeta_{\chi;P;\pi_1,\ldots,\pi_n}(s)=\exp\sum\limits_{d=1}^\infty\frac{q^{-ds}}d\Phi_{q^d}(\chi^{(d)};P;\pi_1^{(d)},\ldots,\pi_n^{(d)}),
\end{equation}
where the characters $\chi^{(d)}:\F_{q^d}\to\C^\times$ and $\pi_1^{(d)},\ldots,\pi_n^{(d)}:\F_{q^d}^\times\to\C^\times$ are defined as 
\begin{equation}
\chi^{(d)}(t)=\chi(\tr t),\ \ \pi_j^{(d)}(t)=\pi_j(Nt),\ \ j=1,\ldots,n,\ \  t\in\F_{q^d};
\end{equation}
here $\tr:\F_{q^d}\to\F_q$ is the trace map and $N:\F_{q^d}\to\F_q$ is the norm map for the extension of fields $\F_q\subset\F_{q^d}$.
Then for $P$ given by (5) and for $\pi_1,\ldots,\pi_{k+n}$ given by (6), equalities (7, 8) imply that the hypergeometric zeta function equals to a ratio 
of products of zeta functions of algebraic varieties with shifted arguments. 

The hypergeometric zeta function can be also represented by a product, 
\begin{equation}
\zeta_{\chi;P;\pi_1,\ldots,\pi_n}(s)=\prod\limits_{x\in\T_{\F_q}^n}(1-\rho(x)\#k(x)^{-s})^{-1},
\end{equation} 
where 
\begin{equation}
\T_{\F_q}^n=(\G_m)_{\F_q}^n=\bigcup\limits_{d=1}^\infty(\F_{q^d}^\times)^n/\Gal(\F_{q^d}/\F_q)
\end{equation}
 is the $n$-dimensional torus over $\F_q$, and for a point 
\begin{equation}
x=[(t_1,\ldots,t_n)\mod\Gal(\F_{q^{d(x)}}/\F_q)]\in\T^n_{\F_q},\ t_1,\ldots,t_n\in\F_{q^{d(x)}}^\times, 
\end{equation}
we have
\begin{equation}
\rho(x)=\chi(\tr P(t_1,\ldots,t_n))\pi_1(Nt_1)\ldots\pi_n(Nt_n).
\end{equation}

{\it The problem is}: to what extent the Weil conjectures can be generalized for hypergeometric zeta functions?

\section{Noncommutative schemes}

Here we provide a new approach to noncommutative algebraic geometry. 

Algebraic varieties can be defined in two ways. The first way is through systems of algebraic equations (algebraic formulas), and the second way is through the points functor
(the functor of sets of points of an algebraic variety with values in various fields and commutative rings).

When we study the problem of generalizing algebraic geometry to noncommutative case, we can start with examples of noncommutative algebraic varieties. We consider the following examples
of noncommutative points functors.

1. The first functor assigns to a skew field $K$ the group $\GL(n,K)$ of invertible $n\times n$-matrices with entries from $K$. We call this example by the group of invertible $n\times n$-matrices
with noncommutative entries, and denote it by $\GL(n)$. 

2. The second functor assigns to a skew field $K$ the set $\P_L^n(K)$ of 1-dimensional left vector subspaces in the standard $(n+1)$-dimensional left vector space $K^{n+1}$ over $K$. We call it by the left $n$-dimensional 
noncommutative projective space, and denote it by $\P_L^n$. Similarly we define the right noncommutative projective space, and denote it by $\P_R^n$.

3. The left (right) noncommutative Grassmannian $\Gr_L^{n,k}$ (resp. $\Gr_R^{n,k}$) assigns to a skew field $K$ the set $\Gr_L^{n,k}(K)$ (resp. $\Gr_R^{n,k}(K)$) 
of $k$-dimensional left (resp. right) vector subspaces in the standard $n$-dimensional left (resp. right) vector space $K^n$ over $K$.

4. Similarly for $0<k_1<\ldots<k_l<n$ one defines the left (right) noncommutative flag variety $\Fl_L^{n; k_1,\ldots,k_l}$ (resp. $\Fl_R^{n;k_1,\ldots,k_l}$) parameterizing flags $0\subset V_1\subset\ldots\subset V_l\subset K^n$
of left (resp. right) vector subspaces $V_j$, $\dim V_j=k_j$, $j=1,\ldots,l$, in the standard $n$-dimensional left (resp. right) vector space $K^n$. 
\medskip

Let us try to describe these varieties by algebraic formulas.

{\it Example.} The inverse matrix to a $2\times 2$-matrix $X=\left(\begin{array}{cc}a_1\ b_1\\ a_2\ b_2\end{array}\right)$ with noncommutative entries is given by the formula
\begin{equation}
X^{-1}=\left(\begin{array}{cc}(a_1-b_1b_2^{-1}a_2)^{-1} \  (a_2-b_2b_1^{-1}a_1)^{-1} \\  (b_1-a_1a_2^{-1}b_2)^{-1} \  (b_2-a_2a_1^{-1}b_1)^{-1} \end{array}\right).
\end{equation}

Formally this formula, for $a_1,b_1,a_2,b_2$ from a skew field $K$, holds only for $a_1,b_1,a_2,b_2\ne 0$ and for $a_1-b_1b_2^{-1}a_2, a_2-b_2b_1^{-1}a_1,b_1-a_1a_2^{-1}b_2,
b_2-a_2a_1^{-1}b_1\ne0$. However, for example, if $a_2, b_2 \ne 0$ and $a_1-b_1b_2^{-1}a_2, b_1-a_1a_2^{-1}b_2\ne 0$, then we have the equality
\begin{equation}
(a_1-b_1b_2^{-1}a_2)^{-1}=-a_2^{-1}b_2(b_1-a_1a_2^{-1}b_2)^{-1}.
\end{equation}
The left hand side of this equality is defined for 
\begin{equation}
b_2\ne 0,\ \  a_1-b_1b_2^{-1}a_2\ne 0,
\end{equation}
 while the right hand side is defined for 
\begin{equation}
a_2\ne 0,\ \  b_1-a_1a_2^{-1}b_2\ne 0.
\end{equation}

This means that we must construct algebraic functions, like (19), on the variety of invertible $2\times 2$-matrices with noncommutative entries by glueing together the objects defined for, say, (20) and for (21).

Similar picture holds for other examples of noncommutative varieties listed above. The local algebraic formulas for algebraic functions on these varieties are given by theory of  
quasideterminants due to I.~M.~Gelfand, V.~S.~Retakh et al. [8].
Quasideterminants of a square matrix with noncommutative entries are noncommutative rational functions of matrix entries. 
The main difficulty in studying noncommutative rational functions is that they contain a lot of nested inversions, and representation of a rational
function by a formula containing nested inversions is not unique. The theory of quasideterminants gives a lot of deep identities involving 
various representations of noncommutative rational functions by formulas with nested inversions. This theory yields a rich algebraic and geometric structure which can be called by 
a noncommutative generalization of classical algebra and geometry (the epoch before Riemann, Hilbert, and Poincar\'e). The question is to 
extend this theory to a noncommutative generalization of modern (commutative) algebra and geometry, 
i.~e. to a noncommutative algebraic geometry which is an abstract framework for the theory of quasideterminants.

{\it The problem is}: to construct a category of objects named {\it noncommutative schemes}, which contains:

(i) the full subcategory equivalent 
to the category of commutative Grothendieck schemes;

(ii) the full subcategory of objects $\Spec A$ called {\it spectra} of noncommutative rings $A$, which is equivalent to the opposite category of the category of noncommutative rings;

(iii) objects $S$ denoted by $\GL(n)$, $\P_L^n$, $\P_R^n$, $\Gr_L^{n,k}$,
$\Gr_R^{n,k}$, $\Fl_L^{n;k_1,\ldots,k_l}$, $\Fl_R^{n;k_1,\ldots,k_l}$ such that for skew fields $K$ the points functors, $S(K)=$ the set of morphisms $\Spec K\to S$,
are equivalent to Examples 1--4 above. 

It seems that computations in [8] concerning noncommutative Vi\'ete theorem and noncommutative symmetric functions can be given a simple interpretation in terms of 
theory of noncommutative schemes and related noncommutative representations of groups.

\end{document}